\documentclass[11]{article}

\usepackage{amsfonts, amsmath,amssymb,latexsym}

\newtheorem{thm}{Theorem}[section]

\def\be#1\ee{\begin{equation}#1\end{equation}}

\newcommand{\bea}{\begin{eqnarray}}
\newcommand{\eea}{\end{eqnarray}}
\newcommand{\beaa}{\begin{eqnarray*}}
\newcommand{\eeaa}{\end{eqnarray*}}
\newcommand{\bei}{\begin{itemize}}
\newcommand{\eei}{\end{itemize}}

\newcommand{\al}{\alpha}
\newcommand{\e}{\varepsilon}
\newcommand{\kk}{K}
\newcommand{\p}{\varphi}
\newcommand{\s}{\sigma}
\newcommand{\ta}{\theta}
\newcommand{\z}{\zeta}
\newcommand{\la}{\lambda}
\newcommand{\La}{\Lambda}
\newcommand{\bk}{\mathbf{k}}
\newcommand{\Z}{\mathbb{Z}}
\newcommand{\X}{\mathbb{X}}
\newcommand{\PP}{\mathbb{P}}
\newcommand{\R}{\mathbb{R}}
\newcommand{\N}{\mathbb{N}}

\newcommand{\E}{\mathbb{E}}
\newcommand{\EE}{\mathcal{E}}
\newcommand{\I}{\mathbb{I}}
\newcommand{\KK}{\mathcal{K}}
\begin{document}
\title{Dependence on the Dimension for Complexity of~Approximation~of~Random~Fields
\footnote{This work is supported by the RFBR under grant number 05-01-00911 and RFBR-DFG under grant number 04-01-04000.}}
\author{ Nora Serdyukova  \footnote{Weierstrass-Institut f\"ur Angewandte Analysis und Stochastik, Mohrenstr. 39, 10117 Berlin, Germany. E-mail:~serdyuko@wias-berlin.de } }

\renewcommand{\baselinestretch}{2}
\bigskip

\author{ Nora Serdyukova  \footnote{Weierstrass-Institut f\"ur Angewandte Analysis und Stochastik, Mohrenstr. 39, 10117 Berlin, Germany. E-mail:~serdyuko@wias-berlin.de } }

\date{}
\maketitle

\begin{abstract}
We consider the $\e$-approximation by $n$-term partial sums of the Karhunen-Lo\`eve expansion to $d$-parametric random fields of tensor product-type in the average case setting. We investigate the behavior, as $d\to \infty$, of the information complexity $n(\e,d)$ of approximation  with error not exceeding a given level $\e$. It was recently shown by M.~A.~Lifshits and E.~V.~Tulyakova that for this problem one observes the curse of dimensionality (intractability) phenomenon. The aim of this paper is to give the exact asymptotic expression
for the information complexity $n(\e,d)$.

\end{abstract}

\noindent

{\bf Key words:}\ random fields, Gaussian
processes, linear approximation error, information-based complexity, tractability, curse of
dimensionality, multivariate linear problems, Karhunen-Lo\`eve expansion.
\bigskip\bigskip
\baselineskip=6.0mm

\section{Introduction}

\par Suppose we have a random function $X(t)$, with $t$ in some compact
parametric set $T$, admitting a series representation via random
variables $\xi_k$ and the deterministic real functions $\p_k$,
namely, 
\[ X(t) = \sum_{k=1}^{\infty}\xi_k \p_k(t),\]
where the series converges in the mean and a.s.\ for each $t \in T$.
A more precise description will be given later.
For any finite set of positive integers $\kk\subset \N$  let $
X_\kk(t) = \sum_{k\in \kk} \xi_k \p_k(t)$. In many problems one
needs to approximate $X$, for instance under $L_2$-norm, with finite rank process $X_\kk$. Natural
questions arise then: how large should be $\kk$ that yields a
given small approximation error? Given the size of $\kk$, which
$\kk$ provides the smallest error? 
\par In this article we address the
first of these questions for a specific class of random functions,
namely {\it tensor product-type random fields} with
high-dimensional parameter sets. The tensor product-type field is 
a separable zero-mean random function $X=\{ X(t)\}_{t\in T}$, $T \subset \R^d$ 
with rectangular parametric set $T$ and covariance function $\KK^{(d)}$
which can be decomposed in a product of equal ``marginal'' covariances depending on different arguments.
Namely, let $T=[0,1]^d$  and 
\be \label{K} \KK^{(d)}(s,t) = \prod_{l=1}^d \KK_l(s_l,t_l) \ee 
for all $s_l, t_l \in [0,1]$,
$s=(s_1, ... , s_d)$, $t = (t_1, ... , t_d)$. Obviously, the integral operator
with the kernel~\eqref{K} is the tensor product of the integral operators with the kernels
 $\KK_l(s_l,t_l)$.
\medskip

Let $\{\lambda_{i}\}_{i\ge 1}$ be a non-negative sequence satisfying
\be  \label{l2}
  \sum_{i=1}^\infty \lambda_{i}^2 <\infty
\ee
and let $\{\varphi_i\}_{i>0}$ be an orthonormal system in $L_2[0,1]$.

Consider a family of tensor product-type random fields \be 
\X=\left\{ X^{(d)}(t), t \in[0,1]^d \right\}\,,\;\; d=1,2,\ldots\, .\ee 
According to the multiparametric Karhunen-Lo\`eve expansion (see~\cite{A} for details),
the following equality in distribution holds
\bea \label{xd1} X^{(d)}(t)
    &=& \sum_{\bk\in \N^d} \xi_k \prod_{l=1}^d \la_{k_l} \prod_{l=1}^d \varphi_{k_l}(t_l)\\
&=&  \sum_{k_1=1}^{\infty} \cdots \sum_{k_d=1}^{\infty}
\xi_{k_1,\ldots,k_d} \la_{k_1} \cdots \la_{k_d} \varphi_{k_1}(t_1)
\cdots \varphi_{k_d}(t_d) ,\;\; \nonumber \eea

\noindent where the series converges a.s.\ for every $t = (t_1, \ldots , t_d)\in [0,1]^d$. 
The collection $\{\xi_k\}$ is an array of non-correlated random variables with zero mean and
unit variance and $\la_{k_l}^2$ and $\varphi_{k_l}$ are, respectively, the eigenvalues and eigenfunctions of the family of integral equations
\[ \la_{k_l}^2 \varphi_{k_l}(t_l) = \int_0^1 \KK_l(s_l,t_l)\varphi_{k_l}(s_l) \mathrm{d}s_l \;,\;\;\;\;\; t_l \in [0,1] \;,\;\;\;\;\; l=1,...,d, \] corresponding to the ``marginal'' covariance operators. Obviously, under assumption~(\ref{l2}) the sample paths of
$X^{(d)}$ belong to $L_2([0,1]^d)$ almost surely and the covariance
operator of $X^{(d)}$ has the system of eigenvalues
\be\label{sobd}
\lambda^2_{\bk}:= \prod_{l=1}^{d} \la_{k_l}^2\,,\,\,\,\,  \bk \in \N^d.
\ee
\par As it was mentioned in~\cite{Sab}, the Karhunen-Lo\`eve expansion or the proper orthogonal decomposition of random functions was introduced independently and almost simultaneously by D.~D.~Kosambi~\cite{Kos}, M.~Lo\`eve~\cite{L}, K.~Karhunen~\cite{Karh1} and~\cite{Karh2},
A.~M.~Obukhov~\cite{Ob} and V.~S.~Pougachev~\cite{Pug}.

\par In the following we drop the index $d$ and write $X(t)$
instead of $X^{(d)}(t)$. 
For any $n>0$, let $X_n$ be the partial sum of (\ref{xd1})
corresponding to $n$ maximal eigenvalues. We study the {\it average case error }
of approximation to $X$ by $X_n$ $$
       e(X, X_n; d) = \left(\E ||X-X_n ||^2_{L_2(T)}
       \right)^{1/2},
$$ as $d \to \infty$. Since in the following we consider only
$L_2(T)$-norms, we will write $|| \cdot||$ instead of
$||\cdot ||_{L_2(T)}$. It is well known (see, for example,
\cite{BS}, \cite{KL} or \cite{R}) that $X_n$ provides the minimal
average quadratic error among all linear approximations to $X$
having rank~$n$.

As we are going to explore a {\it family} of random functions, it is
more natural to investigate {\it relative} errors, that is to compare
the error size with the size of the function itself.

Let
\[ \La:=\sum_{i=1}^\infty \la_{i}^2, 
\]
then
\[
\E\| X\|^2 = \sum_{\bk\in\N^d} \la_{\bk}^2 = \La^d.
\]

Then the {\it average case information complexity} for the normalized error criterion reads as the minimal number of terms in $X_n$ (or, equivalently, of maximal eigenvalues, if they would be ordered) needed to approximate $X$ with the error not exceeding a given level $\e$:
\[
n(\e,d) :=
 \min \{n : \frac{e(X, X_n; d)}{\left(\E\| X\|^2\right)^{1/2}} \leq \e\}
= \min \{n : \E \|X-X_n \|^2 \leq \e^2\La^d \}.
\]

The study of $n(\e,d)$  we are interested in here belongs to the
class of problems dealing with the dependence of the information complexity
for linear multivariate problems on the dimension, see the works
of H.~Wo\'{z}niakovski (\cite{W92}, \cite{W94a}, \cite{W94b},
\cite{W2006}) and the references therein.


It was suggested in \cite{LT} to use an auxiliary probabilistic
construction for studying the properties of deterministic array of
eigenvalues (\ref{sobd}). We follow this approach.

Consider a sequence  of independent identically distributed random
variables $\left\{U_l\right\}, \,\,l=1,2,...$ with the common
distribution given by
\be\label{def_Ul}
\PP(U_l=-\log
\la_{i})=\frac{\la_{i}^2}{\La} \,,\,\,\,\,i=1,2,...
\ee

Under the assumption
\be\label{3d_moment}
\sum_{i=1}^{\infty}|\log \la_{i}|^3 \la_{i}^2 \;<\; \infty,
\ee
the condition $ \E | U_l |^3 < \infty$ is obviously satisfied.

Let $M$ and $\s^2$ denote, respectively, the mean and the variance of $U_l$. Clearly,
\bea \nonumber M &=& - \sum_{i=1}^{\infty}\log
\la_{i}\,\frac{ \la_{i}^2}{\La},
\\ \nonumber
\s^2 &=&
\sum_{i=1}^{\infty}|\log \la_{i}|^2\, \frac{ \la_{i}^2}{\La}\; -\;
M^2.
\eea
Then the third central moment of $U_l$ is given by
\[
\al^3 := \E(U_l -M)^3 =  - \sum_{i=1}^{\infty}\left(\log
\la_{i}\right)^3\,\frac{ \la_{i}^2}{\La} \; - \; 3M\s^2 \;-\;M^3.
\]
If \eqref{3d_moment} is verified, we have $|M|<\infty$, $0 \leq \s^2
< \infty$ and $|\al|< \infty.$

In the sequel the explosion coefficient \be\label{Expl} \EE := \La
e^{2M} \ee will play a significant role, because its contribution into the ``curse of dimensionality'' is the largest. It was shown in \cite{LT}
that by concavity of the logarithmic function $\EE >1$, except for the totally degenerate case when the
number of strictly positive eigenvalues is zero or one. In other
words, $\EE =1$ iff $\s=0$. Henceforth we will exclude this
degenerate case.

The following result was obtained in \cite{LT}, Theorem 3.2.

\begin{thm}\label{T:LT_3.2}
Assume that the sequence $\{\la_{i}\}_{i\geq1}$  satisfies the condition
\[  \sum_{i=1}^{\infty}|\log \la_{i}|^2\,
\la_{i}^2 < \infty.
\]
Then for every $\e \in (0,1)$ we have
\[
\lim_{d \to \infty} \frac{\log n(\e,d) - d \log \EE}{\sqrt{d}} =
2q,
\]
where the quantile $q=q(\e)$ is chosen
from the equation
\be\label{q}
1 - \Phi\left( \frac{q}{\s}\right) = \e^2.
\ee
\end{thm}

The authors of \cite{LT} conjectured that under further assumptions on
the sequence $\{\la_{i}\}$ one can prove that
\[
n(\e,d) \approx \frac{C(\e) \EE^d e^{2q\sqrt{d}}}{\sqrt{d}}\ ,\;
\;\, d \to \infty.
\]
We are going to confirm this conjecture.

\section{Main result}

It turns out that two different cases depending on the nature
of the distribution of $U_l$ should be distinguished. The proof and the final result depend on whether  this
distribution  is a {\it lattice} one or not.

Recall that one calls a discrete distribution of a
random variable $U$ a lattice distribution, if there exist numbers $a$
and $h>0$ such that every possible value of $U$ can be
represented in the form $a+\nu h$, where $\nu$ is an integer.
The number $h$ is called a span of the distribution. In the following, when studying the
lattice case, we assume that $h$ is a maximal span of the
distribution, i.e.\ one cannot represent all possible values of
 $U_l$ in the form $b+\nu h_1$ for some $b$ and $h_1>h$.

Definition \eqref{def_Ul} yields that  the variables $U_l$ have a common lattice
distribution iff $\la_{i} = C e^{-n_i h}$ for some positive $C$, $h$
and $n_i \in \N$. We call this situation the {\it lattice case}
and will assume that $h$ is chosen to be the largest possible.
Otherwise we say that the {\it non-lattice case} takes place.

\par By $f(d)=o(g(d))$ we mean that $\lim_{d \to \infty}\frac{f(d)}{g(d)} =0$. In particular, $ f(d)=g(d)\left( 1 + o(1) \right) $ means that $\lim_{d \to \infty}\frac{f(d)}{g(d)} =1 $.

\bigskip
\begin{thm}
Let the sequence $ \left\{ \la_{i}\right\}_{i \geq 1}$
satisfy~\eqref{3d_moment}. \par Then for every $\e \in (0,1)$ it holds 
\[ n(\e,d) = K\ \phi(\frac{q}{\s})\  \EE^d e^{2q \sqrt{d}}
\,d^{-1/2} \left( 1 + o(1) \right),\; \;\, d \to \infty,\] where
\beaa\phi(x) &=& \frac{1}{\sqrt{2 \pi}} \, e^{-x^2/2 },\\
K &=& \begin{cases}
\frac{h}{\s(1-e^{-2h})} & \mathrm{in \;the\; lattice \;case,} \\
 \frac{1}{2\s}  &\mathrm{otherwise},
\end{cases}
\eeaa and the quantile $q=q(\e)$ is defined in~\eqref{q}.
\end{thm}

\begin{remark}
\bei \item One can see that the complexity of approximation
increases exponentially as $d \to \infty$. This phenomenon is
referred to as the {\it curse of
dimensionality} or { \it intractability}, see
e.g.~\cite{R}~and~\cite{W94a}. The notion of ``curse of
dimensionality'' dates back at least to Bellman~\cite{B}.
\item By l'Hospital's rule
\[ \lim_{h \to 0} \frac{h}{\s \left(1-e^{-2h}\right)} = \frac{1}{2\s }\ . \]
\eei
\bigskip
\end{remark}

\begin{pf}

Let $\z = \z(\e, d)$ be the maximal positive number such that
the sum of eigenvalues satisfies
\[
\sum_{\bk\in\N^d : \la_{\bk}< \z}\la_{\bk}^2   \leq \e^2\La^d.
\]
Define a lattice set in $\N^d$
\[
\mathrm{A} =  \mathrm{A}(\e, d):= \left\{  \bk \in \N^d : \la_{\bk} \geq \z
\right\} = \left\{  \bk \in \N^d : \prod_{l=1}^{d} \la_{k_l} \geq \z \right\}.
\]
Since for any $\bk \in \mathrm{A}$ it holds that $\la_{\bk} >0$, one can
write \beaa \nonumber && n(\e,d) \; = \; card
\left(\mathrm{A}\right) \;
 =\;
\sum_{\bk \in \mathrm{A}}\frac{\la_{\bk}^2}{\la_{\bk}^2}
\\ \nonumber
&=& \sum_{ \bk \in
\N^d : -\sum \log \la_{k_l} \leq -\log \z}\La^d \exp\Big\{-2
\sum_{l=1}^d \log \la_{k_l}\Big\}\prod_{l=1}^{d}\PP(U_l=-\log
\la_{k_l})
\\ \nonumber
&=& \La^d \,\E\exp \Big\{2 \sum_{l=1}^d U_l  \Big\}
\I_{\{\sum_{l=1}^d U_l \leq -\log \z \}}.
\eeaa
For  centered and normalized sums
\[ Z_d = \frac{\sum_{l=1}^d U_l
- d M }{\s\sqrt{d}}
\]
we have
\[ \left\{\sum_{l=1}^d U_l \leq -\log
\z \right\} = \left\{ Z_d \leq \theta \right\},
\]
where
\be\label{ta_def}
\ta = \ta(\e,d) = - \frac{\log\z + d M}{\s \sqrt{d}}.
\ee

We show now that $\ta$ has a useful probabilistic meaning in terms of $\{U_l\}$ and of their sums.
Applying Lemma~3.1 from \cite{LT} we have for any $d \in \N$ and $z \in \R^1$
 \beaa
 \sum_{\bk\in \N^d :\la_{\bk} < z} \la_{\bk} ^2 \; & = & \La^d\;\PP \left( \sum_{l=1}^d U_l > -\log z
 \right)\\ =  \La^d\;\PP \left( Z_d > - \frac {\log z + d M}{\s
 \sqrt{d}}\right) &=&  \La^d\;\PP \left( Z_d > \ta_z   \right),
 \eeaa
 where
$$ \ta_z = -\frac{\log z + d M}{\s \sqrt{d}}.
$$
Fix $\e \in (0,1)$. Observe that
$$   \sum_{\bk\in \N^d :\la_{\bk} < z} \la_{\bk} ^2 \; \leq \; \e^2 \La^d
$$
iff
$$  \PP \left( Z_d > \ta_z   \right) \; \leq \; \e^2 .
$$
Therefore,
$\ta =\ta(\e, d)$ defined by \eqref{ta_def} is the $(1-\e^2)$-quantile of the
distribution of $Z_d$, namely,
\[
\ta(\e, d)=\min\{\ta:\ \PP\left(Z_d>\ta\right)\le \e^2\}
 =\min\{\ta:\ \PP\left(Z_d\le \ta\right) >1- \e^2\}.
\]
Let $q=q(\e)$ be the quantile of the normal distribution function chosen from
the equation~\eqref{q}.  Then in view of the Central Limit Theorem
\be\label{ta}
\ta(d,\e) \to \frac{q(\e)}{\s}\; ,
\;\; d \to \infty,
\ee
for any fixed $\e \in(0,1)$.


Now let us return to the information complexity. We
obtain \bea \nonumber n(\e,d) &=& \EE^d \, \E \exp \{ 2 \s
\sqrt{d} Z_d\} \I_{\{Z_d \leq \ta \}}
\\ \nonumber
&=& \EE^d \, \exp \{ 2 \s \sqrt{d}\ta \}
\int_{-\infty}^{\ta} \exp \{ 2 \s \sqrt{d}(z-\ta) \} \,\mathrm{d}F_d(z) ,
\eea
where $F_d(z) = \PP(Z_d < z)$ and $\EE$ is defined in~\eqref{Expl}.

Denote
\[
 \Psi_d(z):= \exp \{ 2 \s \sqrt{d}(z-\ta) \}
\]
and integrate by parts the integral
\[
 \int_{-\infty}^{\ta} \Psi_d(z)  \,\mathrm{d}[F_d(z)-F_d(\ta)] =
 \int_{-\infty}^{\ta} [- F_d(z)+F_d(\ta)] \,\mathrm{d}\Psi_d(z).
\]

From now on we have to distinguish the lattice and non-lattice cases.
\vspace{0.3cm}

\begin{bf}  Non-lattice case \end{bf}
\par In the following part of the proof we will assume that the
distribution of $\left\{ U_l \right\}$ is not a lattice one. This is
true in the most interesting cases, such as the Brownian sheet (the Wiener-Chentsov random field),  the completely tucked Brownian sheet (the Brownian pillow), the d-variate Hoeffding, Blum, Kiefer and Rosenblatt process (see Appendix for details).

\par In view of~\eqref{3d_moment} we are able to apply the Cram\'{e}r-Esseen Theorem (cf.\ Theorem~2
\S 42 in~\cite{GK},  Theorem~5.21 \S 5.7 of Chapter V
in~\cite{Pe2} or Theorem~4 \S 3 of Chapter VI
 in~\cite{Pe1}). It leads to
\bea \label{CE_th} \nonumber \int_{-\infty}^{\ta} [-
F_d(z)+F_d(\ta)] \,\mathrm{d}\Psi_d(z)& =& \int_{-\infty}^{\ta} [-
\Phi(z)+\Phi(\ta)] \,\mathrm{d}\Psi_d(z)
\\
+ \frac{\al^3}{6\s^3
\sqrt{2\pi d}}\int_{-\infty}^{\ta} [(z^2-1)e^{-z^2/2} &-&
((\ta^2-1)e^{-\ta^2/2} ] \,\mathrm{d}\Psi_d(z) + o\left(
\frac{1}{\sqrt{d}}\right)
\\ \nonumber
= I_1 + I_2 -I_3  - I_4 &+& o\left( \frac{1}{\sqrt{d}}\right),
\eea where \beaa \nonumber I_1 &=& \int_{-\infty}^{\ta} [-
\Phi(z)+\Phi(\ta)] \,\mathrm{d}\Psi_d(z),
\\
I_2 &=& \frac{\al^3}{6\s^3 \sqrt{2\pi
d}}\int_{-\infty}^{\ta}  z^2 e^{-z^2/2}  \,\mathrm{d}\Psi_d(z),
\\ \nonumber
I_3 &=& \frac{\al^3}{6\s^3 \sqrt{2\pi d}}\int_{-\infty}^{\ta}
e^{-z^2/2} \,\mathrm{d}\Psi_d(z),
\\ \nonumber
I_4 &=& \frac{\al^3}{6\s^3
\sqrt{2\pi d}}\left(\ta^2 - 1 \right) e^{-\ta^2/2} =
\frac{\al^3}{6\s^3 \sqrt{2\pi d}}\left(\left(\frac{q}{\s}\right)^2
- 1 \right) \exp\left\{  -\frac{ q^2 }{ 2\s^2 } \right\}\left( 1 + o(1) \right).
\eeaa
The last equivalence is provided by~\eqref{ta}.

Since $\mathrm{d}\Psi_d(z) = 2 \s \sqrt{d}\Psi_d(z)d z$,  the integral $I_2$ is
given, after a change of variable, by
\beaa \nonumber
 I_2 & = & I_2(d,\ta) =
 \frac{\al^3}{3\s^2 \sqrt{2
 \pi d }}\int_0^{\infty}(\ta - \frac{y}{\sqrt{d}})^2 \exp\{
 -\frac{1}{2}(\ta - \frac{y}{\sqrt{d}})^2\}\exp\{-2\s y \} \; \mathrm{d} y
 \eeaa
 where $y = -\sqrt{d}(z-\ta)$.

\par  For any $ d =1,2,...$

\[
0\le
\left(\ta - \frac{y}{ \sqrt{d}}\right)^2 \exp\{
-\frac{1}{2}(\ta - \frac{y}{ \sqrt{d}})^2 \}
\leq (|\ta| + y)^2.
\]
This estimate gives us the majorant required in the Lebesgue's dominated convergence theorem. Using~\eqref{ta} and passing to the limit in the integral
we obtain, as $d \to \infty$,
\[
I_2(d,\ta) =
 \frac{\al^3}{6\s^3 \sqrt{2
 \pi d }} \;\left(\frac{q}{\s}\right)^2 \exp\left\{  -\frac{ q^2 }{ 2\s^2 } \right\} \left( 1 +
 o(1)\right).
\]
Similarly,
\[
 I_3(d,\ta) =
 \frac{\al^3}{6\s^3 \sqrt{2
 \pi d }}  \exp\left\{  -\frac{ q^2 }{ 2\s^2 } \right\} \left( 1 + o(1)\right).
\]
Thus we obtain that $ \sqrt{d}I_4 =  \sqrt{d}(I_2 -I_3 )\left( 1 + o(1)\right)$,
hence, $I_2-I_3-I_4 =  o\left( \frac{1}{\sqrt{d}}\right)$.

Consider the main integral $I_1$.

\bea \label{I1} \nonumber I_1 &=& I_1(d, \ta) = \int_{-\infty}^\ta
[-\Phi(z)+\Phi(\ta)] \,\mathrm{d}\Psi_d(z)\\
\nonumber &=&\frac{1}{\sqrt{2 \pi}} \int_{-\infty}^\ta  \exp\{ 2 \s \sqrt{d}(z-\ta) \} \exp\{ -z^2/2\}  \; \mathrm{d} z \\
\nonumber &=&\frac{1}{\sqrt{2 \pi d }}\int_0^{\infty}\exp\{
 -\frac{1}{2}(\ta - \frac{y}{\sqrt{d}})^2\}\exp\{-2\s y \} \; \mathrm{d} y
\\ &=&  \frac{1}{2\s \sqrt{2 \pi d}} \exp\left\{  -\frac{ q^2 }{ 2\s^2 } \right\}\left( 1 + o(1) \right)\;,\;\; d \to \infty. \eea

 Then
\[
n(\e,d) = \frac{\EE^d\, \exp\{ 2q \sqrt{d} \}}{2 \s \sqrt{d}}
\frac{1}{\sqrt{2 \pi}} \exp\left\{  -\frac{ q^2 }{ 2\s^2 } \right\} \left( 1 + o(1)
\right),
\]
as asserted. \vspace{0.3cm}

\begin{bf} Lattice case \end{bf}
\par Now we will proceed under the assumption that the random variables $U_l$
have a lattice distribution. Let possible values of the random variable $U_l$ be
\[
\tilde{a} + \nu h , \; \nu = 0,\pm 1, \pm 2,...
\] where
$\tilde{a} = M+a$ is a shift, and  $h$ is the maximal span of
the distribution. Therefore, all possible values of $Z_d$ have the
form
\[
\frac{d a + \nu h}{\s \sqrt{d}}, \; \nu = 0,\pm 1, \pm 2,...
\]
Introduce the function \[ S(x) = [x] - x+ \frac{1}{2}, \] where
$[x]$ denotes, as usual, the integer part of $x$, and consider
\[
S_d(x) = \frac{h}{\s} \,S\left( \frac{x \s \sqrt{d} - d
a}{h}\right).
\]
Let $F_d(z)$ be as above. Then under the assumption \eqref{3d_moment}
Esseen's result (see Theorem~1~\S~43 in~\cite{GK}) yields
\[
F_d(z) - \Phi(z) = \frac{e^{-z^2/2}}{\sqrt{2\pi}} \left(
\frac{S_d(z)}{\sqrt{d}} - \frac{\al^3 (z^2 - 1)}{6 \s^3
\sqrt{d}}\right) + o \left( \frac{1}{\sqrt{d}}\right)
\]
uniformly in $z$.

Comparing with \eqref{CE_th}, we observe that one only needs to
evaluate the additional term \beaa J &=&\frac{1}{\sqrt{2 \pi d}}
\int_{-\infty}^{\ta}[- S_d(z) e^{-z^2/2} + S_d(\ta) e^{-\ta^2/2}]
\mathrm{d} \Psi_d(z)
\\
&=& \frac{1}{\sqrt{2 \pi d}} \int_{-\infty}^{\ta} \Psi_d(z)
\mathrm{d} \left( S_d(z) e^{-z^2/2}\right) =  J_1 - J_2 + J_3 ,
\eeaa where \beaa \nonumber J_1 &=& \frac{1}{\sqrt{2 \pi d}}
\int_{-\infty}^{\ta} \Psi_d(z) S_d'(z)e^{-z^2/2} \mathrm{d} z,
\\ \nonumber J_2 &=& \frac{1}{\sqrt{2 \pi d}}
\int_{-\infty}^{\ta} \Psi_d(z) S_d(z) z e^{-z^2/2} \mathrm{d} z,
\eeaa  and $J_3$ is a ``discrete part'', which  is defined in the
following way. Notice that $S(x)$ is a periodic function with period one,
therefore $S_d(x)$ possesses the period $h/\s \sqrt{d}$ and has
jumps at points $\{ \frac{kh + da}{\s \sqrt{d}}, k \in \Z \}$. If
the point $\ta$ belongs to this lattice then there exists an integer $k'$
such that $\ta = \frac{k'h + da}{\s \sqrt{d}}$. Hence, one can
integrate the discontinuous  part of the integral $J$ with respect
to the measure $\frac{h}{\s}\delta_{\frac{kh + da}{\s \sqrt{d}}}$
and obtain
$$ J_3 = \frac{1}{\sqrt{2 \pi d}}\, \frac{h}{\s} \sum_{k=-\infty}^{k'}
\Psi_d\left(\frac{kh + da}{\s \sqrt{d}}\right)
\exp\{-\frac{1}{2}\left(\frac{kh + da}{\s \sqrt{d}}\right)^2\}.$$

\par We start with the estimation of $J_1$. At the points where the
derivative $S_d'(z)$ makes sense, one can easy calculate that
$S_d'(z) = \frac{h}{\s}S\left(\frac{z\s\sqrt{d} - da}{h} \right) =
-\sqrt{d} $, therefore, similarly to the non-lattice case, by
the Lebesgue's dominated convergence theorem we have \bea \label{J1}\nonumber J_1 &=&
\frac{-\sqrt{d}}{\sqrt{2 \pi d}} \int_{-\infty}^{\ta} \exp\{ 2\s
\sqrt{d}(z-\ta)\} \exp\{ -z^2/2\} \mathrm{d} z \\\nonumber &=&
\frac{-1}{\sqrt{2 \pi d}} \int_{0}^{\infty} \exp
\{-\frac{1}{2}(\ta
- \frac{y}{\sqrt{d}})^2 \}\exp\{-2\s y \} \mathrm{d} y\\
&=& \frac{-1}{ 2\s \sqrt{2 \pi d} }\exp\left\{  -\frac{ q^2 }{ 2\s^2 } \right\}\left( 1 + o(1) \right)\;,\; d \to
\infty , \eea and it yields $\sqrt{d}J_1 = - \sqrt{d} I_1\left( 1+
o(1)\right)$.

As for the integral $J_2$, this one, as $d$ is large enough, becomes negligible. 
Indeed, \beaa J_2 &=& \frac{1}{\sqrt{2 \pi d}}
\int_{-\infty}^{\ta}
\exp\{ 2\s \sqrt{d}(z-\ta)\}S_d(z) z  \exp\{ -z^2/2\} \mathrm{d} z \\
 &=&       \frac{1}{\sqrt{2 \pi d}}  \frac{1}{\sqrt{d}} \int_{0}^{\infty} \exp \{-\frac{1}{2}(\ta
- \frac{y}{\sqrt{d}})^2 \} (\ta - \frac{y}{\sqrt{d}})S_d(\ta
- \frac{y}{\sqrt{d}})\exp\{-2\s y \} \mathrm{d}y \\
&\leq& \frac{3 h}{2 \s d\sqrt{2 \pi } }\int_{0}^{\infty} \exp
\{-\frac{1}{2}(\ta - \frac{y}{\sqrt{d}})^2 \} (\ta -
\frac{y}{\sqrt{d}})\exp\{-2\s y \} \mathrm{d}y
\\ &=& \frac{3h}{ 4\s^2 d \sqrt{2 \pi }  }\left( \frac{q}{\s}\right)^2 \exp\left\{  -\frac{ q^2 }{ 2\s^2 } \right\}\left( 1 + o(1) \right)\;,\; d \to \infty
. \eeaa And, of course,  $J_2 = o\left(
\frac{1}{\sqrt{d}}\right)$.

\par Now we consider the most
essential summand \bea \label{J3} \nonumber J_3 &=&
\frac{1}{\sqrt{2 \pi d}}\, \frac{h}{\s} \sum_{k=-\infty}^{k'}
\exp\{2\s \sqrt{d} \left( \frac{k h + d a }{\s \sqrt{d}} - \ta
\right) \} \exp\{- \frac{1}{2} \left(
 \frac{k h + d a}{\s \sqrt{d}} \right)^2\}
\\\nonumber&=& \frac{1}{\sqrt{2 \pi d}}\,
\frac{h}{\s} \sum_{k=-\infty}^{k'}\exp\{2 h (k - k') \} \exp\{-
\frac{1}{2} \left(
 \frac{k h + d a}{\s \sqrt{d}} \right)^2\}
 \\ \nonumber&=& \frac{1}{\sqrt{2 \pi d}}\,
\frac{h}{\s} \sum_{l=0}^{\infty}  \exp\{-2 h l\} \exp\{-
\frac{1}{2} \left(
 \frac{(k' - l) h + d a}{\s \sqrt{d}} \right)^2\}
 \\\nonumber &=& \frac{1}{\sqrt{2 \pi d}}\,
\frac{h}{\s} \sum_{l=0}^{\infty}  \exp\{-2 h l\} \exp\{-
\frac{1}{2} \left( \ta - \frac{l h}{\s \sqrt{d}}\right)^2
\\ &=& \frac{1}{\s
\sqrt{ d}}\, \frac{h }{(1-e^{-2h})} \frac{1}{\sqrt{2
\pi}}\exp\left\{  -\frac{ q^2 }{ 2\s^2 } \right\}\left( 1 + o(1) \right)\;,\; d \to \infty. \eea We obtain
$$\sqrt{d}J_3 = \sqrt{d} \frac{2h }{(1-e^{-2h})} I_1 \left( 1 + o(1)\right).$$

\par Putting together \eqref{I1}, \eqref{J1} and \eqref{J3}, we
get  \[ n(\e,d) = \frac{\EE^d \, e^{ 2q \sqrt{d} }}{ \s
\sqrt{d}} \frac{h }{(1-e^{-2h})} \frac{1}{\sqrt{2 \pi}} \,
\exp \left\{  -\frac{ q^2 }{ 2\s^2 } \right\} \left( 1 + o(1) \right),\; \;\, d \to \infty. \]

\end{pf}
\ $\Box$

\section{Appendix. Examples of tensor product-type \\random fields}
This section contains some examples of random fields to which the above general result can be applied.
\subsection{Wiener-Chentsov random field}
The Wiener-Chentsov field (the Brownian sheet) (see \cite{Lft}) is a zero-mean Gaussian random function $W^{(d)}$ with the covariance function equal to a product of the covariance functions corresponding to the Wiener process  $W$:
\[ \KK_{W^{(d)}} (s,t) = \prod_{l=1}^d \min\{s_l, t_l\} ,\;s=(s_1,...,s_d),\;t=(t_1,...,t_d) \in T.   \]
Therefore the marginal eigenvalues have the following form:
 \[\la_{W;i}^2 = (\pi (i-1/2))^{-2},\;i=1,2,\ldots \,.\]

\subsection{Completely tucked Brownian sheet}
The completely tucked Brownian sheet (the Brownian pillow) is a zero-mean Gaussian random function $B^{(2)}$ with the covariance function, equal to a product of the covariance functions corresponding to the standard Brownian bridge $B(t) = W(t) - t W(1)$, namely
\[ \KK_{B^{(2)}} (s,t) = \prod_{l=1}^2\left( \min\{s_l, t_l\} - s_l t_l  \right) ,\;s,t \in [0,1]^2.   \]
Correspondingly, the marginal eigenvalues (see.~\cite{AnDar}) are equal to \[\la_{B;i}^2 = (\pi
i)^{-2},\;i=1,2,\ldots \,.\]
\par  In the literature different terms are in use for this random field. In~\cite{vW} the term ``completely
tucked Brownian sheet'' is used; in~\cite{CsHor} -- ``tied-down Kiefer
 process''; in~\cite{KonProt} this field is called ``the Brownian pillow''.

The notion of ``completely tucked Brownian sheet'' and its generalization for the case $d>2$,
 was introduced by J. R.~Blum, J.~Kiefer and M.~Rosenblatt~\cite{BKRos}
as the limiting distribution for a functional of empirical process occurring in nonparametric testing of independency, so-called ``the independence
empirical process''~(see~\cite{vW}). Therefore the $d$-parametric generalization of the completely tucked Brownian sheet is often referred to as ``the~d-variate Hoeffding, Blum, Kiefer and Rosenblatt process'' (see, for example,~\cite{KonProt}).
The mention of Hoeffding's name in the term is motivated by the fact that the test studied in~\cite{BKRos}
was equivalent to the one suggested earlier by W.~Hoeffding in~\cite{Hoef}. But the limiting distribution, the covariance function, the eigenvalues and the eigenfunctions of the corresponding integral equation were obtained in~\cite{BKRos}. Higher-dimensional generalizations were later treated in~\cite{Dugue} and in~\cite{Deh}.

\subsection{Centered Gaussian processes}
In some statistical problems it is convenient to use centered empirical processes and corresponding limiting Gaussian processes.

\par For any Gaussian process $X = \{ X(t)\}$, $t
\in[0,1]$ we define the centered process
\[ \mathring{X}(t) := X(t)-\int_0^1 X(u) \mathrm{d}u.  \]
The centered Brownian bridge $\mathring{B}$, also referred to in the literature as the Watson process, was introduced in~\cite{Wats} for nonparametric goodness-of-fit testing on a circle. G.~S.~Watson showed that the covariance function is given by
\[ \KK_{\mathring{B}}(s,t) =  \min\{s,t\} - st + \frac{1}{2}(s^2 + t^2-s-t)+\frac{1}{12}     \;,\;\;s,\,t \in[0,1],\]
and the covariance operator with this kernel has a double spectrum, i.e.
\[ \la_{\mathring{B};2i}^2 = \la_{\mathring{B};(2i-1)}^2 =(2\pi i)^{-2},\;i=1,2,\ldots \, .\]

\par The covariance function of the centered Wiener process $\mathring{W}$ has the form
\[ \KK_{\mathring{W}}(s,t) =  \min\{s,t\}  +\frac{1}{2}(s^2 +t^2)-s-t +\frac{1}{3}   \;,\;\;s,\,t \in[0,1],\]
and the corresponding eigenvalues coincide with those of the standard Brownian bridge, i.e.
\[ \la_{\mathring{W};i}^2  = \la_{B;i}^2 = (\pi i)^{-2},\;i=1,2,\ldots, \]
that is in accordance with the well-known equality in distribution for $L_2$-norms of the Brownian bridge
and centered Wiener process, see~\cite{BegNOrs}.

\par Centered integrated Brownian bridge
\[ \breve{B}(t) = \Bar{B}(t) - \int_0^1\Bar{B}(u)\mathrm{d}u,\] where
\[ \Bar{B}(t) = \int_0^t B(u)\mathrm{d}u ,\,t \in[0,1]\] was considered in a framework of goodness-of-fit testing and small deviation probabilities under $L_2$-norm in~\cite{HenzeN} and in~\cite{BegNOrs}, where
its covariance function

\[ \KK_{\breve{B}}(s,t) =  \frac{st\min\{s,t\}}{2} - \frac{\min\{s,t\}^3}{6} -
       \frac{(st)^2}{4} - \frac{s^2+t^2}{6} - \frac{s^4+t^4}{24} + \frac{s^3 +t^3}{6}
       + \frac{1}{45}, \]
 \noindent $s,\,t \in[0,1] $ and eigenvalues \[ \la_{\breve{B};i}^2 = (\pi i)^{-4}
 ,\;i=1,2,\ldots\,
 \] 
\noindent were obtained.

\subsection{Multivariate extensions of the Anderson-Darling process}

\par The tensor product of Anderson-Darling processes $A^{(d)}(t)$, $t~\in~[0,1]^d$
is a zero-mean Gaussian random function $A^{(d)}(t)$, $t~\in~[0,1]^d$ with the covariance function
\[ \KK_{A^{(d)}}(s,t) = \prod_{l=1}^d
    \frac{\min\{s_l, t_l\} - s_l t_l}{ \sqrt{s_l (1 - s_l)} \sqrt{t_l (1 - t_l)} }
        ,\;s_l, t_l \in [0,1].   \]
\noindent The eigenvalues of the corresponding covariance operator are given by
\[ \la_{\bk}^2 = \prod_{l=1}^d \frac{1}{k_l (k_l +1) } , \; \bk=(k_1, \ldots , k_d)  \in \N^d. \]

 \par In the one-dimensional case the Anderson-Darling process coincides in distribution with $\frac{B(t)}{\sqrt{t(1-t)}}$ , $t \in[0,1]$ and 
was introduced in~\cite{AnDar} in the context of goodness-of-fit testing. T.~Anderson and D.~Darling obtained its covariance function and the exact spectrum.

\par In~\cite{Pycke} another multivariate extension of Anderson-Darling process, defined as a zero-mean Gaussian process with the covariance function
\[ \KK_A^{\mu}(s,t) = \left( \frac{\min\{s, t\} - s t}{ \sqrt{s (1 - s)} \sqrt{t (1 - t)}
}\right)^{\mu}
        ,\;s,\, t \in [0,1] ,\, \mu >0,     \]is given.

\noindent The eigenvalues of its covariance operator are of the form
\[ \la_{\mu,j}^{2} =\frac{\mu}{(\mu +j-1)(\mu + j)}\,, \; j=1,2, \ldots \,.  \]
When parameter $\mu$ is positive integer, the random field, defined in such a way, (more precisely, the square of its $L_2$-norm) is the limiting distribution
for Cram\'er-von Mises type statistics.

\bigskip
\bigskip

\noindent {\Large\bf Acknowledgements} \bigskip

\noindent
This paper was partially written while the author was visiting
the Institut f\"{u}r Matematische
Stochastik, Georg-August-Universit\"at, G\"ottingen. Special
thanks are due to Professor M.A. Lifshits for formulation of the problem
and constant
encouragement and to Professor M. Denker for his support and for
providing excellent working conditions.

\end{document}